\documentclass[11pt]{article}
\usepackage{latexsym}
\usepackage[all]{xy}
\usepackage{amsmath}
\usepackage{amssymb}
\usepackage{amsfonts}

\newtheorem{thm}{Th\'eor\`eme}[section]
\newtheorem{prop}[thm]{Proposition}
\newtheorem{lem}[thm]{Lemme}

\newtheorem{df}[thm]{D\'efinition}
\newtheorem{cor}[thm]{Corollaire}
\newtheorem{conj}[thm]{Conjecture}

\begin{document}

\title{\textbf{D\'enombrabilit\'e des classes d'\'equivalences d\'eriv\'ees de vari\'et\'es alg\'ebriques}}
\bigskip
\bigskip

\author{
Mathieu Anel \\
\small{Department of Mathematics, Middlesex College}\\
\small{The University of Western Ontario} \\
\small{London, Ontario N6A 5B7} \\
\small{Canada}
\bigskip\\
Bertrand To\"en\\
\small{Laboratoire Emile Picard, UMR CNRS 5580}\\
\small{Universit\'{e} Paul Sabatier, Bat 1R2} \\
\small{Toulouse Cedex 09}\\
\small{France}}

\date{July 2007}

\maketitle

\begin{abstract}
Soient $S$ un sch\'ema affine,
$X \longrightarrow S$ une famille miniverselle de sch\'emas projectifs et lisses, et 
$D$ une cat\'egorie triangul\'ee fix\'ee. On d\'emontre que les points
$s\in S$ tels que la cat\'egorie d\'eriv\'ee de la fibre en $s$, $D_{coh}^{b}(X_{s})$, soit \'equivalente \`a $D$,
forment un ensemble au plus d\'enombrable. Nous d\'eduisons de cela que l'ensemble des classes
d'isomorphisme
des vari\'et\'es complexes lisses et projectives qui poss\`edent une cat\'egorie
d\'eriv\'ee fix\'ee est au plus d\'enombrable.
Notre d\'emonstration
passe par la construction d'un certain pr\'echamp classifiant les dg-cat\'egories satur\'ees et connexes, ainsi qu'une
application des p\'eriodes allant du champ des vari\'et\'es lisses et projectives vers ce pr\'echamp, et 
qui \`a une vari\'et\'e associe un dg-mod\`ele pour sa cat\'egorie d\'eriv\'ee. 
\end{abstract}

\medskip

\tableofcontents

\bigskip

\section*{Introduction} 

\bigskip

\`A toute vari\'et\'e alg\'ebrique projective et lisse $X$ (disons sur un corps $k$) 
on peut associer $D_{coh}^{b}(X)$, sa cat\'egorie
d\'eriv\'ee coh\'erente born\'ee. Pour deux telles vari\'et\'es $X$ et $Y$ \`a fibr\'es
canoniques (ou anti-canoniques) amples, on sait que $X$ et $Y$ sont isomorphes si et seulement si
les cat\'egories $D_{coh}^{b}(X)$ et $D_{coh}^{b}(Y)$ sont \'equivalentes en tant que
cat\'egories triangul\'ees (voir \cite{bo}). Cependant, il existe en g\'en\'eral des
exemples o\`u $X$ et $Y$ ne sont pas isomorphes mais o\`u les cat\'egories 
d\'eriv\'ees $D_{coh}^{b}(X)$ et $D_{coh}^{b}(Y)$ sont \'equivalentes (voir \cite{ro} pour 
des r\'ef\'erences). 
Cela pose \'evidemment la question de la classification des vari\'et\'es projectives et lisses
$X$ qui poss\`edent une cat\'egorie d\'eriv\'ee fix\'ee (\`a \'equivalence triangul\'ee pr\`es). 
Dans cette direction Y. Kawamata propose la conjecture suivante (voir \cite{ro} pour une discussion 
de cette conjecture). 

\begin{conj}\label{conj}
Soit $D$ une cat\'egorie triangul\'ee $k$-lin\'eaire. Alors, l'ensemble
des classes d'iso\-morphisme de vari\'et\'es projectives et lisses $X$ sur $k$ 
telles que $D_{coh}^{b}(X)$ soit \'equivalente \`a $D$ (en tant que
cat\'egorie $k$-lin\'eaire triangul\'ee) est fini. 
\end{conj}

L'objet de ce travail est de proposer une approche g\'eom\'etrique de cette conjecture qui 
utilise un certain \emph{espace de modules de cat\'egories triangul\'ees}. L'id\'ee g\'en\'erale, et na\"{\i}ve, 
est de chercher \`a construire un tel espace de modules $\mathcal{M}$ qui param\'etrise les
cat\'egories triangul\'ees, ainsi qu'une certaine \emph{application 
des p\'eriodes}
$$\Pi : \mathcal{V} \longrightarrow \mathcal{M},$$
o\`u $\mathcal{V}$ est un espace de modules pour les vari\'et\'es projectives et lisses, et o\`u 
$\Pi$ envoie $X$ sur $D_{coh}^{b}(X)$. On cherchera alors \`a montrer que $\Pi$ 
est infinit\'esimalement injective (i.e. non-ramifi\'ee), et donc \`a fibres discr\`etes. Dans le cas o\`u ces 
espaces de modules sont assez proches d'\^etre des vari\'et\'es alg\'ebriques ses fibres
seraient donc finies. Cela impliquerait \ref{conj} car la fibre de $\Pi$ pris au point
$D$ est pr\'ecis\'ement l'ensemble dont la conjecture \ref{conj} pr\'edit la finitude. \\

Dans ce travail nous construisons des mod\`eles aux espaces $\mathcal{V}$ et $\mathcal{M}$
et au morphisme $\Pi$, et nous montrons que $\Pi$ est non-ramifi\'e en un certain sens. Avant de passer
aux d\'etails techniques des d\'efinitions de ces objets signalons qu'une cons\'equence 
de leurs existences est le th\'eor\`eme suivant, qui est l'\'enonc\'e principal de cet article. 

\begin{thm}{(Cor. \ref{c2})}\label{ti}
Soit $X \longrightarrow S$ une famille miniverselle de sch\'emas projectifs, lisses et 
g\'eom\'etriquement connexes, avec $S$ un sch\'ema de type fini sur un corps $k$. 
Soit $D$ une cat\'egorie triangul\'ee $k$-lin\'eaire, et 
$S(D)$ le sous-ensemble de $S(k)$ form\'e des 
points $s$ tels que $D_{coh}^{b}(X_{s})$ soit
\'equivalente \`a $D$ (comme cat\'egorie triangul\'ee $k$-lin\'eaire).
Alors $S(D)$ est un ensemble au plus d\'enombrable. 
\end{thm}

Il faut noter que la conjecture \ref{conj} implique l'\'enonc\'e pr\'ec\'edent (voir notre
remarque au \S 8). \\

Bien que \ref{ti} est relativement loin de la conjecture \ref{conj}, 
 il affirme qu'elle  est \emph{moralement} vraie si l'on remplace
\emph{fini} par \emph{au plus d\'enombrable} (''moralement'' car 
on ne consid\`ere ici que des vari\'et\'es qui apparaissent dans une m\^eme famille miniverselle). 
De plus, dans le cas des vari\'et\'es alg\'ebriques complexes un argument 
transcendant permet de d\'eduire la d\'enombrabilit\'e des classes
d'\'equivalences d\'eriv\'ees.

\begin{thm}{(Cor. \ref{t2})}\label{ti2}
Soit $D$ une cat\'egorie triangul\'ee $\mathbb{C}$-lin\'eaire. 
Alors, l'ensemble des classes d'isomorphismes
de vari\'et\'es complexes lisses et projectives $X$ telles que
$D_{coh}^{b}(X)$ soit \'equivalente (comme cat\'egorie triangul\'ee
$\mathbb{C}$-lin\'eaire) \`a $D$ est au plus d\'enombrable. 
\end{thm}

Quelques mots sur le contenu de cet article. Le lecteur ne sera pas surpris d'apprendre que nous
avons chercher des mod\`eles pour $\mathcal{V}$, $\mathcal{M}$ et $\Pi$ dans le
cadre des champs. Il le sera peut-\^etre un peu d'apprendre que 
les mod\`eles que nous construisons de ces objets ne sont \emph{pas}
des champs alg\'ebriques, et de plus que le mod\`ele que nous donnons
pour $\mathcal{M}$  est un pr\'echamp qui n'est m\^eme pas un champ. Plus pr\'ecis\'ement, 
pour $\mathcal{V}$ on prendra, comme on peut s'y attendre, le champ 
$\mathcal{VAR}^{c}_{smpr}$ des sch\'emas propres, lisses et g\'eom\'etriquement connexes
(dont les sections au-dessus d'un sch\'ema $S$ est le groupo\"{\i}de des
morphismes $X \longrightarrow S$ propres, lisses et \`a fibres g\'eom\'etriquement connexes). 
Il est bien connu que ce champ n'est pas alg\'ebrique, car il existe des d\'eformations
formelles dans $\mathcal{VAR}^{c}_{smpr}$ qui ne sont pas alg\'ebrisables. 
Comme mod\`ele pour l'objet $\mathcal{M}$ nous proposons le pr\'echamp
$\mathcal{DGC}^{c}_{sat}$ qui param\'etrise les \emph{dg-cat\'egories satur\'ees et connexes, prises
\`a quasi-\'equivalence pr\`es}\footnote{Il est bien connu que la notion
de cat\'egorie triangul\'ee poss\`ede beaucoup de d\'esavantages qui sont r\'esolus 
en la rempla\c{c}ant par la notion de dg-cat\'egorie. C'est seulement au prix de ce remplacement que
l'objet $\mathcal{DGC}^{c}_{sat}$  peut-\^etre d\'efini et poss\`ede une chance de bien 
se comporter. }  (voir \S 2,3 pour plus de d\'etails). Le pr\'echamp 
$\mathcal{DGC}^{c}_{sat}$ n'est pas un champ, et ses pr\'efaisceaux de morphismes
ne sont pas non plus des faisceaux (il s'agit donc d'un pr\'echamp en un sens encore
plus faible que celui de \cite{lm}). Il n'est donc pas alg\'ebrique, et son
champ associ\'e ne l'est probablement pas non plus, et ce pour la m\^eme raison
que pour le cas de $\mathcal{VAR}^{c}_{smpr}$ (\`a savoir l'existence
de d\'eformations formelles non alg\'ebrisables, voir notre remarque 
au \S 8). Enfin, nous construisons 
un morphisme de pr\'echamps
$$\Pi : \mathcal{VAR}^{c}_{smpr} \longrightarrow \mathcal{DGC}^{c}_{sat},$$
qui \`a un sch\'ema propre et lisse $X \longrightarrow S=Spec\, k$ associe
une dg-cat\'egorie $L_{pe}(X)$ de complexes parfaits sur $X$, qui est 
un dg-mod\`ele \`a la cat\'egorie d\'eriv\'ee parfaite $D_{parf}(X)$. 
Le fait que $X \mapsto L_{pe}(X)$ soit compatible aux changements de bases $k\rightarrow k'$, et 
donc le fait que $\Pi$ existe comme morphisme de pr\'echamps, est une
propri\'et\'e sp\'ecifique \`a l'utilisation des dg-cat\'egories qui ne serait pas
v\'erifi\'ee si l'on avait chercher \`a d\'efinir $\mathcal{DGC}^{c}_{sat}$ \`a l'aide
de cat\'egories triangul\'ees. 

Nous montrerons alors que le morphisme $\Pi$ est non-ramifi\'e, au sens o\`u 
une d\'eformation infinit\'esimale dans $\mathcal{VAR}^{c}_{smpr}$ est triviale si et seulement
si son image dans $\mathcal{DGC}^{c}_{sat}$ l'est (voir Prop. \ref{p2}). Cela peut sembler de peu 
d'int\'er\^et \'etant donn\'e que les pr\'echamps $\mathcal{VAR}^{c}_{smpr}$ et 
$\mathcal{DGC}^{c}_{sat}$ ne sont pas des champs alg\'ebriques. Cependant, 
ces deux objets partagent une propri\'et\'e en commun avec les
champs alg\'ebriques qui permet de tirer des cons\'equences
int\'eressantes de l'existence de $\Pi$ et de son caract\`ere non ramifi\'e. 
Cette propri\'et\'e est la \emph{quasi-repr\'esentabilit\'e}, ou en d'autres termes
la repr\'esentabilit\'e du morphisme diagonal. Pour $\mathcal{VAR}^{c}_{smpr}$
la quasi-repr\'esentabilit\'e est une cons\'equence de l'existence
des espaces alg\'ebriques de Hilbert, et est donc bien connue. La quasi-repr\'esentabilit\'e
de $\mathcal{DGC}^{c}_{sat}$ (tout au moins \`a un proc\'ed\'e de faisceautisation pr\`es) 
est elle une cons\'equence des r\'esultats de \cite{tova} qui affirment l'existence d'un
champ alg\'ebrique classifiant les objets d'une dg-cat\'egorie satur\'ee (voir Prop. \ref{p1}). On montre aussi 
que la diagonale de $\mathcal{DGC}^{c}_{sat}$ poss\`ede une propri\'et\'e 
de \emph{quasi-compacit\'e d\'enombrable}, qui affirme que 
l'intersection de deux affines au-dessus de $\mathcal{DGC}^{c}_{sat}$ poss\`ede un 
recouvrement par un nombre d\'enombrable de sch\'emas affines (voir Prop. \ref{p1}).
Le th\'eor\`eme \ref{ti} est alors une cons\'equence formelle du caract\`ere
non ramifi\'e de $\Pi$, de la quasi-repr\'esentabilit\'e de $\mathcal{DGC}^{c}_{sat}$, 
de la propri\'et\'e de quasi-compacit\'e d\'enombrable, et d'un th\'eor\`eme
de D. Orlov affirmant que deux vari\'et\'es lisses et projectives 
poss\`edent des cat\'egories d\'eriv\'ees \'equivalentes si et seulement
leur dg-mod\`eles sont quasi-\'equivalents (voir \cite{or}). Au passage, nous d\'emontrons une version
plus g\'en\'erale du th\'eor\`eme \ref{ti}, valable pour les familles propres et lisses, mais pour
laquelle la notion de cat\'egorie triangul\'ee doit \^etre remplac\'ee par celle
de dg-cat\'egorie.  \\

\bigskip

\textbf{Remerciements:} Nous remercions B. Keller et G. Vezzosi pour plusieurs
discussion sur la th\'eorie des d\'eformations des dg-cat\'egories qui ont influenc\'e ce travail (bien que
cette th\'eorie des d\'eformation ait pu \^etre soigneusement \'evit\'ee). 
Un grand merci \`a C. Simpson pour nous avoir expliqu\'e comment 
le th\'eor\`eme \ref{ti2} d\'ecoulait du th\'eor\`eme \ref{ti}. 
Nous remercions enfin J. Lurie pour sa remarque sur la non repr\'esentabilit\'e
du champ des dg-cat\'egories satur\'ees, et C. Voisin 
pour ses commentaires. \\

\bigskip

\textbf{Conventions et notations:} Nous notons $CAlg$ la cat\'egorie des anneaux commutatifs.
Pour $k\in CAlg$, nous notons $k-CAlg$ celle
des $k$-alg\`ebres commutatives. Lorsque nous consid\'ererons des faisceaux ou bien 
des champs, la cat\'egorie $Aff$ des sch\'emas affines sera
munie de la topologie \'etale. Le produit fibr\'e homotopique de pr\'echamps, aussi appel\'e
\emph{2-produit fibr\'e},  sera not\'e $-\times^{h}_{-}-$. Enfin, on renvoie 
\`a \cite{ke,to1,to2,tova} pour des rappels sur le "monde d\'eg\'e" (dg-alg\`ebres, dg-cat\'egories, 
dg-modules, dg-modules parfaits \dots). \\

\bigskip

\section{Dg-cat\'egories satur\'ees et connexes}

Dans ce paragraphe nous rappellerons bri\`evement les notions
de dg-cat\'egories propres, lisses, triangul\'ees et satur\'ees. On renvoie \`a
\cite{ke,tova} pour des d\'efinitions plus d\'etaill\'ees. \\

Fixons un anneau commutatif de base $k$.
Une dg-cat\'egorie (sur $k$) est une cat\'egorie enrichie dans la
cat\'egorie mono\"{\i}dale des complexes (non born\'es) de $k$-modules.
\`A une telle dg-cat\'egorie $T$ on associe une cat\'egorie
homotopique $[T]$, dont les objets sont les m\^emes que ceux 
de $T$ et dont les ensembles de morphismes sont les $H^{0}$
des complexes de morphismes de $T$. Il existe une notion naturelle
de morphisme de dg-cat\'egories qui est celle de foncteur enrichi. 
Un tel foncteur, $T \longrightarrow T'$, est une \emph{quasi-\'equivalence} 
s'il induit des quasi-isomorphismes
au niveau des complexes de morphismes et s'il induit un foncteur
essentiellement surjectif $[T] \longrightarrow [T']$. On s'int\'eressera
particuli\`erement \`a la cat\'egorie homotopique des dg-cat\'egories (au-dessus de $k$), 
not\'ee $Ho(k-dg-cat)$, obtenu en localisant celle des dg-cat\'egories le long
des quasi-\'equivalences (voir \cite{tab} pour une justification de l'existence raisonable d'une telle localisation).

Pour une dg-cat\'egorie $T$, on dispose d'une cat\'egorie
de $T^{op}$-dg-modules $T^{op}-Mod$, et de sa cat\'egorie homotopique
$D(T^{op}):=Ho(T^{op}-Mod)$. La cat\'egorie $D(T^{op})$ est naturellement munie
d'une structure triangul\'ee.
 Le plongement de Yoneda enrichi permet de
construire un foncteur pleinement fid\`ele
$$[T] \longrightarrow D(T^{op}).$$
Nous dirons alors que $T$ est \emph{triangul\'ee} si l'image essentielle de ce foncteur
consiste en la sous-cat\'egorie $D(T^{op})_{c}$ des objets compacts 
dans $D(T^{op})$. 

Nous dirons qu'une dg-cat\'egorie $T$ est \emph{compactement engendr\'ee} s'il existe une dg-alg\`ebre $B$ sur $k$, 
et une \'equivalence de cat\'egories triangul\'ee $ D(T^{op})\simeq D(B)$.  De mani\`ere
\'equivalente, $T$ est compactement engendr\'ee si la cat\'egorie triangul\'ee $D(T^{op})$ poss\`ede
un g\'en\'erateur compact. Nous dirons que $T$ est \emph{localement propre} si pour toute
paire d'objets $(x,y)$ dans $T$, le complexe $T(x,y)$ est un complexe parfait de $k$-modules.
Nous dirons alors que $T$ est \emph{propre} si elle est compactement engendr\'ee et localement propre. 
De mani\`ere \'equivalente $T$ est propre si elle est Morita \'equivalente (au sens d\'eriv\'ee)
\`a une dg-alg\`ebre $B$ dont le complexe sous-jacent est un complexe parfait de $k$-modules. 

Pour deux dg-cat\'egories $T$ et $T'$ on peut d\'efinir leur produit tensoriel
$T\otimes T'$ (sur $k$), dont l'ensemble d'objets est le produit des ensembles d'objets
de $T$ et $T'$, et dont les complexes de morphismes sont les produits tensoriels
des complexes de morphismes de $T$ et $T'$. Cette construction peut se d\'eriver
\`a gauche pour donner un produit tensoriel d\'eriv\'e $T\otimes^{\mathbb{L}}T'$. 
Ce produit tensoriel d\'eriv\'e munit $Ho(k-dg-cat)$ d'une structure mono\"{\i}dale sym\'etrique
que l'on montre \^etre ferm\'ee (voir \cite{to1}).

Une dg-cat\'egorie $T$ d\'efinit un $T\otimes T^{op}$-module, qui envoie une paire
d'objets $(x,y)$ sur le complexe $T(y,x)$, et donc un objet dans 
$D(T\otimes^{\mathbb{L}}T^{op})$. Nous dirons que $T$ est \emph{lisse}
si l'objet $T$ est compact dans $D(T\otimes^{\mathbb{L}}T^{op})$. \\

Nous pouvons enfin donner la d\'efinition de dg-cat\'egorie satur\'ee.

\begin{df}\label{d1}
Une dg-cat\'egorie $T$ au-dessus de $k$ est \emph{satur\'ee}, si elle
est propre, lisse et triangul\'ee.
\end{df}

Pour tout morphisme d'anneaux commutatifs $k \longrightarrow k'$ il existe 
un foncteur de changement de bases
$$-\otimes_{k}^{\mathbb{L}}k' : Ho(k-dg-cat) \longrightarrow Ho(k'-dg-cat),$$
qui est 
adjoint \`a gauche du foncteur d'oubli. On remarque alors que 
les dg-cat\'egories propres et lisses sont stables par changement de bases. Les dg-cat\'egories
triangul\'ees, et donc satur\'ees, ne le sont pas. Cependant, le foncteur
d'inclusion $Ho(k-dg-cat^{tr}) \hookrightarrow Ho(k-dg-cat)$ 
des dg-cat\'egories triangul\'ees dans les dg-cat\'egories poss\`ede
un adjoint \`a gauche $T \mapsto \widehat{T}_{pe}$ (voir \cite[\S 7]{to1}). On peut donc 
d\'efinir un foncteur de changement de bases pour les dg-cat\'egories 
triangul\'ees par
$$\begin{array}{ccc}
Ho(k-dg-cat^{tr}) & \longrightarrow & Ho(k'-dg-cat^{tr}) \\
T & \longmapsto & \widehat{(T\otimes_{k}^{\mathbb{L}}k')}_{pe}.
\end{array}$$
Ce dernier changement de bases pr\'eserve alors les 
dg-cat\'egories satur\'ees (on utilisera ici les r\'esultats
de \cite[Lem. 2.6]{tova} qui affirment qu'\^etre propre et lisse
est invariant par la construction $T \mapsto \widehat{T}_{pe}$).

Pour terminer signalons que la sous-cat\'egorie pleine
$Ho(k-dg-alg^{sat})\subset Ho(k-dg-cat)$, form\'ee des dg-cat\'egories
satur\'ees, poss\`ede la description suivante. On consid\`ere 
les dg-alg\`ebres $A$ sur $k$ qui sont propres et lisses (i.e.
parfaites comme $k$-dg-modules et comme $A\otimes^{\mathbb{L}}A^{op}$-dg-module). 
Pour deux telles dg-alg\`ebres $A$ et $B$ on note 
$Mor(A,B)$ l'ensemble des classes d'isomorphisme d'objets
compacts dans la cat\'egorie triangul\'ee $D(A\otimes^{\mathbb{L}}B^{op})$. 
Pour trois dg-alg\`ebres propres et lisses $A$, $B$ et $C$ on dispose d'un morphisme de composition
$$Mor(A,B)\times Mor(B,C) \longrightarrow Mor(A,C),$$
qui \`a $M\in D(A\otimes^{\mathbb{L}}B^{op})$ et 
$N\in D(B\otimes^{\mathbb{L}}C^{op})$ associe
$$M\otimes_{B}^{\mathbb{L}}N\in D(A\otimes^{\mathbb{L}}C^{op}).$$
Ceci fait des dg-alg\`ebres propres et lisses une cat\'egorie 
not\'e $Ho(k-dg-alg^{pl}_{mor})$. Les th\'eor\`emes de
\cite{to1} montrent alors qu'il existe une \'equivalence de cat\'egories
$$Ho(k-dg-alg^{pl}_{mor}) \simeq Ho(k-dg-cat^{sat}).$$
Le foncteur qui induit cette \'equivalence envoie 
une dg-alg\`ebre propre et lisse $A$ sur la dg-cat\'egorie
$\widehat{A}_{pe}$, des $A^{op}$-dg-modules compacts.

\section{Le pr\'echamp $\mathcal{DGC}_{sat}^{c}$}

Soient $k$ un anneau commutatif et $T\in Ho(k-dg-cat)$ une dg-cat\'egorie. On d\'efinit la cohomologie de Hochschild
de $T$ 
par la formule suivante (voir \cite[\S 8.1]{to1})
$$HH^{i}_{k}(T):=Ext^{i}_{T\otimes_{k}^{\mathbb{L}}T^{op}}(T,T),$$
o\`u les groupes $Ext^{i}$ sont calcul\'es dans la cat\'egorie triangul\'ee
$D(T\otimes_{k}^{\mathbb{L}}T^{op})$. Comme il se doit $HH_{k}^{*}(T)$ est munie
d'une structure naturelle de $k$-alg\`ebre gradu\'ee. On dispose donc d'un morphisme naturel
$k \longrightarrow HH^{0}(T)$.

Notons au passage que si $A$ est une dg-alg\`ebre et $T:=\widehat{A}_{pe}$ est sa 
dg-cat\'egorie des $A$-dg-modules compacts, alors $HH^{*}(T)$ est naturellement
isomorphe \`a $HH^{*}(A)$ d\'efini \`a l'aide du complexe de Hochschild usuel (voir par exemple
\cite{ke}). 

\begin{df}\label{d2}
Nous dirons qu'une $k$-dg-cat\'egorie $T$ est \emph{connexe} si pour
tout morphisme d'an\-neaux commutatifs $k \longrightarrow k'$ les deux conditions
suivantes sont satisfaites
\begin{enumerate}
\item 
$$HH^{i}_{k'}(T\otimes_{k}^{\mathbb{L}}k')\simeq 0 \quad \forall \; i<0.$$
\item Le morphisme naturel 
$$k' \longrightarrow  HH^{0}_{k'}(T\otimes_{k}^{\mathbb{L}}k')$$
est un isomorphisme. 
\end{enumerate}
\end{df}

Pour tout anneau commutatif $k\in CAlg$ notons 
$\mathcal{DGC}_{sat}^{c}(k)$ le groupo\"{\i}de des 
dg-cat\'egories satur\'ees et connexes dans $Ho(k-dg-cat)$. Pour un morphisme
$k\rightarrow k'$ dans $CAlg$ on a vu que l'on disposait de foncteur
de changement de bases $\widehat{((-)\otimes_{k}^{\mathbb{L}}k')_{pe}}$ qui pr\'eservaient
les dg-cat\'egories satur\'ees, et donc 
de foncteurs naturels
$$\widehat{((-)\otimes_{k}^{\mathbb{L}}k')_{pe}} : \mathcal{DGC}_{sat}^{c}(k) \longrightarrow \mathcal{DGC}_{sat}^{c}(k').$$
Ces changements 
de bases pr\'eservent les dg-cat\'egories satur\'ees et connexes, 
et on a donc d\'efini un pr\'echamp 
$\mathcal{DGC}_{sat}^{c}$
sur le site des sch\'emas affines. 

\begin{df}\label{d3}
Le \emph{pr\'echamp des dg-cat\'egories satur\'ees et connexes} est le 
pr\'echamp $\mathcal{DGC}_{sat}^{c}$ d\'efini ci-dessus. 
\end{df}

Il se trouve que le pr\'echamp $\mathcal{DGC}_{sat}^{c}$ n'est pas un champ (disons pour
la topologie \'etale). Une premi\`ere raison
est que les pr\'efaisceaux d'isomorphismes entre deux objets ne sont pas des faisceaux. En effet, 
ces pr\'efaisceaux sont les pr\'efaisceaux des classes d'isomorphismes d'une $\mathbb{G}_{m}$-gerbe
(voir ci-dessous pour plus de d\'etails) et ne sont donc en g\'en\'eral pas des faisceaux. Cela 
provient en r\'ealit\'e du fait que $\mathcal{DGC}_{sat}^{c}$ a \'et\'e d\'efini comme
le $1$-tronqu\'e d'un $2$-pr\'echamp dont les $1$-pr\'echamps d'isomorphismes sont 
eux des $1$-champs (voir \S 7 et \cite[\S 4.3 (6)]{to3} pour des remarques sur l'aspect champs sup\'erieurs).
Mais de plus, les donn\'ees de descente ne sont pas effectives dans $\mathcal{DGC}_{sat}^{c}$, m\^eme
apr\`es que les pr\'efaisceaux d'isomorphismes aient \'et\'e remplac\'es par des
faisceaux, o\`u m\^eme si $\mathcal{DGC}_{sat}^{c}$ est remplac\'e par
le $2$-pr\'echamp sus-cit\'e. Cela est d\^u au fait qu'il existe des formes tordues non triviales
de dg-cat\'egories pour la topologie \'etale.

\section{Quasi-repr\'esentabilit\'e de $\mathcal{DGC}_{sat}^{c}$}

Nous avons d\'efini un pr\'echamp $\mathcal{DGC}_{sat}^{c}$ que nous avons vu ne pas
\^etre un champ. Il se trouve que le champ associ\'e \`a $\mathcal{DGC}_{sat}^{c}$ n'est 
pas un champ alg\'ebrique (voir notre remarque au \S 8). Nous montrons cependant dans la proposition suivante
que son morphisme diagonal est repr\'esentable.

\begin{prop}\label{p1}
Pour tout sch\'ema affine $X\in Aff$, et toute paire de morphismes de pr\'echamps
$s,t : X \longrightarrow \mathcal{DGC}_{sat}^{c}$, le
faisceau $a(\underline{Iso}(s,t))$, associ\'e au pr\'efaisceau
$$\underline{Iso}(s,t):=\mathcal{DGC}_{sat}^{c}\times_{(\mathcal{DGC}_{sat}^{c}\times \mathcal{DGC}_{sat}^{c})}^{h}X,$$
est repr\'esentable par un espace alg\'ebrique localement de pr\'esentation finie (sur $X$). De plus, 
l'espace alg\'ebrique $a(\underline{Iso}(s,t))$ est une r\'eunion d\'enombrable de sous-espaces
alg\'ebriques ouverts et de pr\'esentation finie (sur $X$). 
\end{prop}

\textit{Preuve:} Soit $X=Spec\, k \in Aff$, et $s,t : X \longrightarrow \mathcal{DGC}_{sat}^{c}$
deux morphismes de pr\'echamps. Les morphismes $s$ et $t$ correspondent \`a deux
objets de $\mathcal{DGC}_{sat}^{c}(k)$, c'est \`a dire \`a deux dg-cat\'egories
$T_{1}$ et $T_{2}$ satur\'ees et connexes. Le pr\'efaisceaux  $\underline{Iso}(s,t)$
se d\'ecrit par le foncteur $k-CAlg \longrightarrow Ens$, qui envoie 
une $k$-alg\`ebre commutative $k'$ sur l'ensemble des
isomorphismes $T_{1}\otimes_{k}^{\mathbb{L}}k' \simeq  T_{2}\otimes_{k}^{\mathbb{L}}k'$
dans $Ho(k'-dg-cat)$.  D'apr\`es les th\'eor\`emes de \cite{to1} le pr\'efaisceau 
$\underline{Iso}(s,t)$ est isomorphe au pr\'efaisceau qui \`a $k'\in k-CAlg$ associe
l'ensembles des classes d'isomorphismes d'objets compacts et inversibles (pour la
composition des morphismes de Morita) dans $D(T_{1}\otimes_{k}^{\mathbb{L}} T_{2}^{op}\otimes_{k}^{\mathbb{L}}k')$.

On pose $T_{0}:=T_{1}\otimes_{k}^{\mathbb{L}}T_{2}^{op}$, qui est une $k$-dg-cat\'egorie propre et lisse 
(mais pas triangul\'ee en g\'en\'eral). 
D'apr\`es \cite{tova}, il existe un ($1$-)champ alg\' ebrique $t_{\leq 0}\mathcal{M}_{T_{0}^{op}}^{simp}$, localement de pr\'esentation finie 
au-dessus de $Spec\, k$, qui
\`a $k'\in A-CAlg$ associe le groupo\"{\i}de des objets compacts $E\in D(T_{0}\otimes_{k}^{\mathbb{L}}k')$ tels que
$Ext^{i}(E,E)=0$ pour $i<0$ et $Ext^{0}(E,E)=k'$. De plus, ce champ alg\'ebrique est une
$\mathbb{G}_{m}$-gerbe au-dessus de son espace de modules $M:=\pi_{0}(t_{\leq 0}\mathcal{M}_{T_{0}^{op}}^{simp})$.
En d'autres termes, l'espace alg\'ebrique $M$ repr\'esente le faisceau associ\'e
au pr\'efaisceau qui \`a $k'\in A-CAlg$ associe l'ensemble des classes d'isomorphismes
d'objets compacts $E\in D(T_{0}\otimes_{k}^{\mathbb{L}}k')$ tels que
$Ext^{i}(E,E)=0$ pour $i<0$ et $Ext^{0}(E,E)=k'$. 

\begin{lem}
Avec les notations ci-dessus, si $E\in D((T_{1}\otimes_{k}^{\mathbb{L}}T_{2}^{op})\otimes_{k}^{\mathbb{L}}k')$ est inversible, 
alors on a 
$$Ext^{i}(E,E)=0 \; \forall \; i<0 \qquad Ext^{0}(E,E)=k'.$$
\end{lem}

\textit{Preuve:} De deux choses l'une: ou bien $T_{1}$ et $T_{2}$ deviennent isomorphes dans le groupo\"{\i}de 
$\mathcal{DGC}_{sat}^{c}(k')$, ou bien elles ne le sont pas. Dans le second cas aucun objet
$E \in D((T_{1}\otimes_{k}^{\mathbb{L}}T_{2}^{op})\otimes_{k}^{\mathbb{L}}k')$ n'est inversible. Supposons alors
que l'on soit dans le premier cas. Nous pouvons donc supposer que $T_{1}=T_{2}=T$. Soit $E \in D((T\otimes_{k}^{\mathbb{L}}T^{op})\otimes_{k}^{\mathbb{L}}k')$ un objet inversible. En multipliant 
par l'inverse de $E$, on trouve une auto-\'equivalence triangul\'ee 
de $D((T\otimes_{k}^{\mathbb{L}}T^{op})\otimes_{k}^{\mathbb{L}}k')$ 
qui envoie $E$ sur $T\otimes_{k}^{\mathbb{L}}k'$.
On a donc
$$Ext^{i}(E,E)\simeq Ext^{i}(T\otimes_{k}^{\mathbb{L}}k',T\otimes_{k}^{\mathbb{L}}k')\simeq HH^{i}(T\otimes_{k}^{\mathbb{L}}k').$$
Ainsi, la condition de connexit\'e sur $T$ implique le lemme. 
\hfill $\Box$ \\

Le lemme pr\'ec\'edent montre que le faisceau $a(\underline{Iso}(s,t))$, associ\'e 
\`a $\underline{Iso}(s,t)$, s'identifie naturellement au sous-faisceau de $M$ form\'e
des objets inversibles. Or, \^etre inversible est une condition ouverte au sens de Zariski (e.g. voir la preuve de \cite[Cor. 3.24]{tova}), et ainsi 
$a(\underline{Iso}(s,t))$ est bien repr\'esentable par un espace alg\'ebrique localement de
pr\'esentation finie sur $Spec\, k$. 

Finalement, on sait que le champ $t_{\leq 0}\mathcal{M}_{T_{0}^{op}}^{simp}$ est une r\'eunion 
d\'enombrable de sous-champs ouverts et de pr\'esentation finie sur 
$Spec\, k$ (voir \cite[\S 3.3]{tova}). Comme l'inclusion
$a(\underline{Iso}(s,t)) \longrightarrow M$ est une 
immersion ouverte 
on en d\'eduit ais\'ement que $a(\underline{Iso}(s,t))$ est encore une r\'eunion d\'enombrable
de sous-espaces alg\'ebriques ouverts et de pr\'esentation finie sur $Spec\, k$.
\hfill $\Box$ \\

\begin{cor}\label{c1}
Pour tout sch\'emas affines $X, Y\in Aff$, et tout morphismes
$$s  : X \longrightarrow \mathcal{DGC}_{sat}^{c} \qquad t : Y \longrightarrow \mathcal{DGC}_{sat}^{c},$$ le
faisceau $a(\underline{Iso}(s,t))$, associ\'e au pr\'efaisceau 
$$\underline{Iso}(s,t):=X\times_{\mathcal{DGC}_{sat}^{c}}^{h}Y,$$
est repr\'esentable par un espace alg\'ebrique localement de type fini (sur $X\times Y$). De plus, 
l'espace alg\'ebrique $a(\underline{Iso}(s,t))$ est une r\'eunion d\'enombrable de sous-espaces
alg\'ebriques ouverts et de pr\'esentation  finie (sur $X\times Y$). 
\end{cor}

\textit{Preuve:} On applique la proposition \ref{p1} au morphisme $X\times Y \longrightarrow \mathcal{DGC}_{sat}^{c}\times \mathcal{DGC}_{sat}^{c}$.
\hfill $\Box$ \\

\section{Application des p\'eriodes}

Notons $\mathcal{VAR}_{smpr}^{c}$ le champ qui \`a $k\in CAlg$ associe 
le groupo\"{\i}de des sch\'emas $X \longrightarrow Spec\, k$, propres, lisses et \`a fibres g\'eom\'etriquement connexes.
Pour $X\longrightarrow Spec\, k$ un objet de $\mathcal{VAR}_{smpr}^{c}$ on peut 
consid\'erer sa dg-cat\'egorie (sur $k$) des complexes parfaits $L_{pe}(X)$, qui est un 
dg-mod\`ele \`a la cat\'egorie triangul\'ee $D_{parf}(X)$ (i.e. il existe une
\'equivalence triangul\'ee naturelle $[L_{pe}(X)]\simeq D_{parf}(X)$, voir \cite[\S 8.3]{to1}). 
On montre de plus que $L_{pe}(X)$ est satur\'ee (voir \cite[Lem. 3.27]{tova}). On montre aussi que
$L_{pe}(X)$ est connexe car on a
$$\widehat{L_{pe}(X)\otimes_{k}^{\mathbb{L}}k'}_{pe}\simeq L_{pe}(X\otimes_{k}k')$$
$$HH^{i}(L_{pe}(X))\simeq Ext^{i}(\Delta_{X},\Delta_{X}),$$
o\`u $\Delta_{X}$ est le faisceau structural de la diagonale et o\`u le
groupe $Ext$ est pris dans $D_{parf}(X\times_{k}X)$.
On d\'efinit ainsi un morphisme de pr\'echamps
$$\Pi : \mathcal{VAR}_{smpr}^{c} \longrightarrow \mathcal{DGC}_{sat}^{c}.$$

\begin{df}\label{d4}
Le morphisme de champs $\Pi$ d\'efini ci-dessus est appel\'e
\emph{application des p\'eriodes}.
\end{df}

\section{Un \'enonc\'e de type Torelli infinit\'esimal}

Soit $K$ un corps alg\'ebriquement clos, $F$ un pr\'echamp et $x\in F(K)$ un point. On rappelle que 
le groupo\"{\i}de tangent de $F$ en $x$, not\'e $T_{x}F$ est d\'efini comme la fibre homotopique du morphisme naturel
$F(K[\epsilon]) \longrightarrow F(K)$ prise au point $x$. En d'autres termes
$$T_{x}F:=F(K[\epsilon])  \times_{F(K)}^{h}\{x\}.$$
L'ensemble des classes d'isomorphismes du groupo\"{\i}de $T_{x}F$ sera appel\'e
\emph{l'espace tangent de $F$ en $x$} et sera not\'e
$$T_{x}^{0}F:=\pi_{0}(T_{x}F).$$
L'espace tangent $T_{x}^{0}F$ est en g\'en\'eral un ensemble point\'e, dont
le point de base sera not\'e $0$ et correspond \`a l'image de $x$ par 
le morphisme naturel $F(K) \longrightarrow F(K[\epsilon])$. \\

La proposition suivante affirme que l'application des p\'eriodes $\Pi$ 
satisfait une propri\'et\'e de type \emph{Torelli infinit\'esimal}.

\begin{prop}\label{p2}
Soit $K$ un corps alg\'ebriquement clos et $X\in \mathcal{VAR}_{smpr}^{c}(K)$ d'image
$\Pi(X)\in \mathcal{DGC}_{sat}^{c}(K)$. Alors le morphisme induit
$$T\Pi : T_{X}^{0}\mathcal{VAR}_{smpr}^{c} \longrightarrow T^{0}_{\Pi(X)}\mathcal{DGC}_{sat}^{c}$$
est tel que 
$$T\Pi^{-1}(0)=\{0\}.$$
\end{prop}

\textit{Preuve:} Soit $u\in T_{X}^{0}\mathcal{VAR}_{smpr}^{c}$, correspondant 
\`a une d\'eformation infinit\'esimale $(X',\alpha)$ de $X$. Dans cette notation
$X' \longrightarrow Spec\, K[\epsilon]$ est un sch\'ema propre et lisse et 
$\alpha : X \longrightarrow X'$ est un morphisme
de $K[\epsilon]$-sch\'emas induisant un isomorphisme
$\alpha_{0} : X\simeq X'\otimes_{K[\epsilon]}K$. Son image par $T\Pi$ est
la d\'eformation infinit\'esimale $(L_{pe}(X'),\gamma(\alpha))$ de $L_{pe}(X)$, o\`u 
$\gamma(\alpha)\in D(X\times_{K[\epsilon]}X')$ est le faisceau structural 
du graphe du morphisme $\alpha$, consid\'er\'e comme un isomorphisme
$$\gamma(\alpha) : L_{pe}(X) \simeq \widehat{(L_{pe}(X')\otimes_{K[\epsilon]}K)_{pe}}$$
dans $Ho(K-dg-cat)$ (voir \cite[Thm. 8.15]{to1}). Supposons que 
$T\Pi(X,\alpha)=0$. En utilisant \cite[Thm. 8.15]{to1} on voit que cela implique qu'il 
existe un objet 
$$E\in D_{parf}(X'\times_{K}X)\simeq D_{parf}(X'\times_{K[\epsilon]}(X\otimes_{K}K[\epsilon])),$$
satisfaisant aux deux conditions suivantes.
\begin{itemize}

\item La transformation de Fourier-Mukai associ\'ee \`a $E$
$$\phi_{E} : D(X') \longrightarrow D(X\otimes_{K}K[\epsilon])$$
est une \'equivalence de cat\'egories.

\item Il existe un isomorphisme dans $D((X'\otimes_{K[\epsilon]}K)\times_{K}X)$
$$E\otimes^{\mathbb{L}}_{K[\epsilon]}K \simeq \gamma(\alpha^{-1}_{0}),$$
o\`u $\gamma(\alpha_{0}^{-1})$ est le faisceau structural du graphe
de l'isomorphisme $\alpha_{0}^{-1}$ (i.e. le transpos\'e de $\gamma(\alpha_{0})$).

\end{itemize}

On commence par remarquer que la seconde condition implique que $E$  est un 
faisceau coh\'erent sur $X'\times_{K}X$, 
plat sur $X'$.  En effet, comme la question est locale sur $X'\times_{K}X$, cela d\'ecoule du lemme suivant.

\begin{lem}\label{lplat}
Soit $A$ un anneau commutatif et $I\subset A$ un id\'eal de
carr\'e nul. Soit $E\in D^{b}(A)$ un complexe born\'e de
$A$-modules. On suppose que le complexe $E\otimes_{A}^{\mathbb{L}}A/I$ est 
cohomologiquement concentr\'e en degr\'e $0$, et de plus que
$H^{0}(E\otimes_{A}^{\mathbb{L}}A/I)$ est un $A/I$-module plat. 
Alors, $E$ est cohomologiquement concentr\'e en degr\'e $0$ et 
de plus $H^{0}(E)$ est un $A$-module plat.
\end{lem}

\textit{Preuve du lemme:} Le lemme de Nakayama implique facilement que l'on a $H^{i}(E)=0$ pour tout $i < 0$.
Montrons que pour tout $A$-module $N$ et tout $i<0$, on a 
$H^{i}(E\otimes^{\mathbb{L}}_{A}N)=0$. Comme tout $A$-module $N$ s'ins\`ere dans une suite
exacte 
$$\xymatrix{ 0 \ar[r] & IN \ar[r] & N \ar[r] & N/IN \ar[r] & 0,}$$ 
on voit qu'il nous suffit de montrer cette assertion pour des $A$-modules
$N$ qui sont des $A/I$-modules. Mais alors, on a
$$E\otimes^{\mathbb{L}}_{A}N \simeq (E\otimes^{\mathbb{L}}_{A}A/I)\otimes_{A/I}^{\mathbb{L}}N.$$
Par hypoth\`ese sur $E$ on trouve donc
$$E\otimes^{\mathbb{L}}_{A}N\simeq (H^{0}(E)\otimes_{A}A/I)\otimes_{A/I}N.$$
Ceci implique bien 
que $H^{i}(E\otimes^{\mathbb{L}}_{A}N)=0$ pour tout $i<0$ et tout $N$.

Lorsque $N=A$ on trouve que $E$ est concentr\'e en degr\'e $0$. De plus, 
pour tout $A$-module $N$ on a alors
$$H^{i}(E\otimes^{\mathbb{L}}_{A}N) \simeq Tor_{-i}^{A}(H^{0}(E),N)=0 \; \forall \; i<0.$$
Ceci implique que $H^{0}(E)$ est plat sur $A$. 
\hfill $\Box$ \\

Comme nous l'avons dit plus haut, le lemme ci-dessus implique que 
$E$ est un faisceau coh\'erent sur $X'\times_{K}X$, 
plat sur $X'$.
De plus, pour tout point $y\in X'$, le 
faisceau coh\'erent $E$ restreint \`a $\{y\}\times_{K}X:=Spec\, k(y) \times_{K} X$
est isomorphe au faisceau gratte-ciel en $(y,\alpha_{0}^{-1}(y))$.
On peut donc repr\'esenter $E$ par un morphisme de $K$-champs
$$E : X' \longrightarrow \underline{Coh}^{1}(X)\simeq X\times_{K} B\mathbb{G}_{m},$$
o\`u $\underline{Coh}^{1}(X)$ est le champ des faisceaux coh\'erents
de longeur $1$ sur $X$. Ce morphisme est tel qu'il existe un isomorphisme naturel
entre le morphisme induit
$$E\otimes_{K[\epsilon]}K : X'\otimes_{K[\epsilon]}K \longrightarrow (X\times_{K} B\mathbb{G}_{m})\otimes_{K[\epsilon]}K\simeq
X\times_{K} B\mathbb{G}_{m},$$
et le morphisme
$$\xymatrix{X'\otimes_{K[\epsilon]}K\ar[r]^-{\alpha_{0}^{-1}} &  X  \ar[r]^-{id\times p} & X \times_{K} B\mathbb{G}_{m},}$$
o\`u $p$ est le $\mathbb{G}_{m}$-torseur trivial sur $X$.
La seconde composante du morphisme $E$ induit un morphisme de $K$-champs 
$$X' \longrightarrow B\mathbb{G}_{m}.$$ 
Ce dernier morphisme d\'efinit alors un fibr\'e en droite sur $X'$ qui est une
d\'eformation infinit\'esimale de $\mathcal{O}_{X'}$, et dont l'image r\'eciproque
sur $X'\times_{K}X$ sera not\'ee $\mathcal{L}$. On consid\`ere alors
$F:=\mathcal{L}^{-1}\otimes E \in D(X'\times_{K}X)$. Cet objet
$F$ est un faisceau coh\'erent sur $X'\times_{K}X$, et par construction
il existe un morphisme de $K$-sch\'emas $f : X'\longrightarrow X$, et un 
isomorphisme de faisceaux coh\'erents $F\simeq \gamma(f)$, o\`u 
$\gamma(f)$ est le graphe du morphisme $f$. Le morphisme
$f$ est alors tel que 
$$f\otimes_{K[\epsilon]}K = \alpha_{0}^{-1} : X'\otimes_{K[\epsilon]}K \simeq X,$$
car ces deux morphismes poss\`edent des graphes dont les faisceaux structuraux
sont isomorphes sur $(X'\otimes_{K[\epsilon]}K)\times_{K}X$.
Ceci montre que 
la d\'eformation $(X',\alpha)$ est triviale. 
\hfill $\Box$ \\

\section{Lieu \`a dg-cat\'egorie fix\'ee des familles miniverselles}

Le th\'eor\`eme qui suit est le r\'esultat principal de ce travail. Rappelons qu'un
morphisme propre et lisse de sch\'emas $X \longrightarrow S$ est une 
famille miniverselle, si pour tout point $s\in S$, le morphisme
de Kodaira-Spencer
$$T_{s}S \longrightarrow H^{1}(X_{s},T_{X_{s}/k(s)})$$
est injectif (comme d'habitude, $X_{s}:=X\times_{S}Spec\, k(s)$ est la fibre
de $X \longrightarrow S$ en $s$).

\begin{thm}\label{t1}
Soit $X \longrightarrow S$ une famille miniverselle de sch\'emas propres, lisses et 
g\'eom\'etri\-quement connexes, avec 
$S$ un sch\'ema quasi-compact. Soit $K$ un corps alg\'ebriquement clos et $T$ une
dg-cat\'egorie sur $K$. Soit $S(T)$ le sous-ensemble de $S(K)$ des points
$Spec\, K \longrightarrow S$ tels que $L_{pe}(X\times_{S}Spec\, K)$ soit
isomorphe dans $Ho(K-dg-cat)$ \`a $T$. Alors $S(T)$ est un ensemble
au plus d\'enombrable.
\end{thm}

\textit{Preuve:} Le morphisme $X \longrightarrow S$ d\'efinit un morphisme de champs
$$S \longrightarrow \mathcal{VAR}_{smpr}^{c},$$
et donc en composant avec $\Pi$ un morphisme de pr\'echamps
$$S \longrightarrow \mathcal{DGC}_{sat}^{c}.$$
La dg-cat\'egorie $T$ d\'efinit un morphisme de pr\'echamps
$$t : Spec\, K \longrightarrow \mathcal{DGC}_{sat}^{c}.$$
Soit $\mathcal{T} \subset \mathcal{DGC}_{sat}^{c}$ le sous-pr\'echamp
plein qui est l'image essentielle (au sens des pr\'echamps) du morphisme $t$. 
Le pr\'echamp $\mathcal{T}$ est naturellement muni d'un morphisme
$\mathcal{T} \longrightarrow Spec\, K$ (i.e. est un $K$-pr\'echamp).
On forme alors le diagramme homotopiquement cart\'esien de pr\'echamps suivant
$$\xymatrix{
S \ar[r]^-{p} & \mathcal{DGC}_{sat}^{c} \\
F \ar[u] \ar[r] & \mathcal{T}. \ar[u]}$$
Soit $Z=a(F)$ le faisceau associ\'e \`a $F$ (pour la topologie \'etale). Le morphisme
$F \longrightarrow S$ se factorise par $Z \longrightarrow S$, et l'image de l'application
$$F(K)\simeq Z(K) \longrightarrow S(K)$$
est par construction le sous-ensemble $S(T)$. 

\begin{lem}\label{l1}
Le faisceau $Z$ est repr\'esentable par un espace alg\'ebrique localement 
de pr\'esentation finie sur $K$. De plus, cet espace alg\'ebrique
est recouvert (au sens de la topologie \'etale) par un nombre 
d\'enombrable de sch\'emas affines.
\end{lem}

\textit{Preuve du lemme \ref{l1}:} On consid\`ere le  diagramme homotopiquement cart\'esien suivant de
pr\'echamps
$$
\xymatrix{
S \ar[r]^-{p} & \mathcal{DGC}_{sat}^{c} \\
F \ar[u] \ar[r] & \mathcal{T} \ar[u] \\
F' \ar[u] \ar[r] & Spec\, K. \ar[u]}$$
Soit $G$ le pr\'efaisceau en groupes des automorphismes du point $t$. On a alors une \'equivalence
de pr\'echamps $\mathcal{T}=BG$. Ainsi, on peut \'ecrire le pr\'efaisceau $F$ comme un quotient
$F'/G$, avec $G$ qui op\`ere sans points fixes sur $F'$. En passant au faisceaux associ\'es, on trouve 
que $Z$ est isomorphe au faisceau quotient de $a(F')$ par $a(G)$, avec $a(G)$ qui op\`ere
toujours librement sur $a(F')$. Or, d'apr\`es la proposition \ref{p1} $a(F')$ est un 
espace alg\'ebrique localement de pr\'esentation et $a(G)$ est un groupe
alg\'ebrique. Ainsi, le faisceau quotient $Z=a(F'/G)$ est bien un 
espace alg\'ebrique localement de pr\'esentation finie. Enfin, comme $a(F')$ 
est une r\'eunion d\'enombrable de sous-espaces ouverts de type fini sur $K$, on voit qu'il
est recouvert par un nombre d\'enombrable de sch\'emas affines. Comme le morphisme
quotient $a(F') \longrightarrow Z$ est un \'epimorphisme, cela implique que
$Z$ est lui-m\^eme recouvert par un nombre d\'enombrable de sch\'emas affines. 
\hfill $\Box$ \\

\begin{lem}\label{l2}
Le $K$-espace alg\'ebrique $Z$ est isomorphe \`a 
un $K$-sch\'ema de la forme $\coprod_{A}Spec\, K$ avec 
$A$ un ensemble au plus d\'enombrable. 
\end{lem}

\textit{Preuve du lemme \ref{l2}:} Soit $x\in F(K)$ un $K$-point. Sa projection $s$ sur
$S$ est telle que $L_{pe}(X_{s})\simeq T$.
L'espace tangent 
de $F$ en $x$ s'inscrit dans un diagramme homotopiquement cart\'esien de groupo\"{\i}des
$$\xymatrix{
T_{s}S \ar[r] & T_{\Pi(X_{s})}\mathcal{DGC}^{c}_{sat} \\
T_{x}F \ar[r] \ar[u] & T_{*}\mathcal{T}.\ar[u] }$$
Comme $\mathcal{T}$ est un sous-pr\'echamp plein de $\mathcal{DGC}^{c}_{sat}$, 
le morphisme $T_{*}\mathcal{T} \longrightarrow T_{T}\mathcal{DGC}^{c}_{sat}$
est pleinement fid\`ele. Ceci montre que $T_{x}F \longrightarrow T_{s}S$ est pleinement fid\`ele, 
et donc que le diagramme induit sur les espaces tangents
$$\xymatrix{
T_{s}S \ar[r]^-{Tp} & T^{0}_{\Pi(X_{s})}\mathcal{DGC}^{c}_{sat} \\
T_{x}F \ar[r] \ar[u] & T^{0}_{*}\mathcal{T}, \ar[u] }$$
est un diagramme cart\'esien d'ensembles point\'es. Or, 
on a clairement $T^{0}_{*}\mathcal{T}\simeq 0$ car $\mathcal{T}$ ne poss\`ede qu'un 
unique objet. Ainsi, $T_{x}F$ s'identifie \`a $Tp^{-1}(0)$. 
Or, le morphisme $Tp$ se factorise en 
$$T_{s}S \longrightarrow T_{X_{s}}^{0}\mathcal{VAR}^{c}_{smpr}\simeq H^{1}(X_{s},T_{X_{s}/k(s)}) \longrightarrow 
T^{0}_{\Pi(X_{s})}\mathcal{DGC}^{c}_{sat}.$$
Le premier de ces morphismes est injectif par hypoth\`ese. Ainsi
$T_{x}F=Tp^{-1}(0)=\{0\}$ d'apr\`es la proposition \ref{p2}. 

Nous venons de voir que $T_{x}F=0$, ce qui implique que 
l'espace alg\'ebrique $Z$ est non-ramifi\'e sur $Spec\ K$ (on remarquera
que $T_{x}F\simeq T_{x}Z$ car $K$ est alg\'ebriquement clos). Ainsi, on trouve
que $Z$ est isomorphe \`a un $K$-sch\'ema de la forme $\coprod_{A}Spec\, K$. 
Enfin, la seconde partie du lemme \ref{l1} implique que $A$ est au plus d\'enombrable.
\hfill $\Box$ \\

Le lemme \ref{l2} implique que $S(T)$, qui est l'image de $Z(K)$ dans $S(K)$,  est un ensemble 
au plus d\'enombrable. 
\hfill $\Box$ \\

\begin{cor}\label{c2}
Soit $X \longrightarrow S$ une famille miniverselle de sch\'emas projectifs, lisses et 
g\'eom\'etri\-quement connexes, avec $S$ un sch\'ema de type fini sur un corps $k$ alg\'ebriquement clos. 
Soit $D$ une cat\'egorie triangul\'ee $k$-lin\'eaire, et 
$S(D)$ le sous-ensemble de $S(k)$ des 
points $s$ tels que $D_{coh}^{b}(X_{s})$ soit
\'equivalente \`a $D$ (comme cat\'egorie triangul\'ee $k$-lin\'eaire).
Alors $S(D)$ est un ensemble au plus d\'enombrable. 
\end{cor}

\textit{Preuve:} D'apr\`es \cite{or}, deux vari\'et\'es lisses et projectives complexes
$X$ et $X'$ sont telles que $D_{coh}^{b}(X)$ et $D^{b}_{coh}(X')$ sont \'equivalentes (comme
cat\'egories triangul\'ees), si et seulement 
si les deux dg-cat\'egories $L_{pe}(X)$ et $L_{pe}(X')$ sont quasi-\'equivalentes. 
\hfill $\Box$ \\

\section{Le cas des vari\'et\'es complexes}

Dans cette section nous d\'emontrerons la cons\'equence suivante du corollaire \ref{c2}. \\

\begin{thm}\label{t2}
Soit $D$ une cat\'egorie triangul\'ee $\mathbb{C}$-lin\'eaire. 
Alors, l'ensemble des classes d'iso\-morphismes
de vari\'et\'es complexes lisses, projectives et connexes $X$ telles que
$D_{coh}^{b}(X)$ soit \'equivalente (comme cat\'egorie triangul\'ee
$\mathbb{C}$-lin\'eaire) \`a $D$ est au plus d\'enombrable. 
\end{thm}

\textit{Preuve:} Il s'agit de construire un nombre d\'enombrables
de familles mini-verselles $\{X^{\alpha} \longrightarrow S^{\alpha}\}$ dont les
fibres donnent toutes les vari\'et\'es complexes lisses, projectives et connexes, puis
d'appliquer le corollaire \ref{c2}. 
Pour cela nous dirons qu'un ensemble de familles
lisses, projectives et connexes $\{X^{\alpha} \longrightarrow S^{\alpha}\}$, avec 
$S^{\alpha}$ des sch\'emas de type fini sur $\mathbb{C}$, est \emph{dominant}
si pour toute vari\'et\'e lisse, projective et connexe $Y$ sur $\mathbb{C}$, il existe
$\alpha$ et $s\in S^{\alpha}(\mathbb{C})$ tels que $Y$ soit 
isomorphe \`a $X_{s}^{\alpha}$. 

\begin{lem}\label{l3}
Il existe un ensemble dominant de familles lisses, projectives et connexes
$\{X^{\alpha} \longrightarrow S^{\alpha}\}_{\alpha \in A}$ avec $A$ un ensemble
d\'enombrable. 
\end{lem}

\textit{Preuve du lemme \ref{l3}:} Soit $H_{n}$ le sch\'ema de Hilbert de $\mathbb{P}^{n}_{\mathbb{C}}$. 
Ce sch\'ema est une r\'eunion disjointe d\'enombrable de vari\'et\'es projectives. 
Soit $H^{li}_{n}$ le sous-sch\'ema ouvert form\'e des sous-sch\'emas
de $\mathbb{P}^{n}_{\mathbb{C}}$ qui sont lisses et connexes sur $\mathbb{C}$. 
Le sch\'ema $H^{li}_{n}$ est lui aussi une r\'eunion disjointe d\'enombrable
de sch\'emas de type fini sur $\mathbb{C}$, que nous noterons
$H^{li}_{n,m}$. De plus, la famille universelle
$$Z_{n,m}\subset H^{li}_{n,m} \times \mathbb{P}_{\mathbb{C}}^{n}$$
induit un morphisme lisse, projectif et connexe $Z_{n,m} \longrightarrow H^{li}_{n,m}$. 
Par construction l'ensemble de familles
$$\{Z_{n,m} \longrightarrow H^{li}_{n,m}\}_{n,m}$$
est dominant et d\'enombrable. \hfill $\Box$ \\

\begin{lem}\label{l4}
Soit $p : X \longrightarrow S$ un morphisme lisse, projectif et connexe avec 
$S$ un sch\'ema de type fini sur $\mathbb{C}$. Alors il existe un
ensemble de familles mini-verselles 
$\{X^{\alpha} \longrightarrow S^{\alpha}\}_{\alpha\in A}$, avec 
$A$ un ensemble d\'enombrable et tel que pour tout $s\in S(\mathbb{C})$
il existe $\alpha$ et $t\in S^{\alpha}(\mathbb{C})$ tels que 
$X_{s}$ et $X^{\alpha}_{t}$ soient isomorphes. 
\end{lem}

\textit{Preuve du lemme \ref{l4}:} En stratifiant $S$ on peut supposer que les propri\'et\'es suivantes
sont satisfaites.

\begin{itemize}

\item Le sch\'ema $S$ est affine, lisse et connexe.

\item Le faisceau coh\'erent $\mathbb{R}^{1}p_{*}(T_{X/S})$ est 
un fibr\'e vectoriel, et sa formation commute aux changements
de bases sur $S$. 

\item Le morphisme de Kodaira-Spencer
$$\Theta : T_{S} \longrightarrow \mathbb{R}^{1}p_{*}(T_{X/S})$$
poss\`ede un noyau $K$ et un conoyau $C$ qui sont des fibr\'es vectoriels sur $S$.
\end{itemize}

Soient $p\in S(\mathbb{C})$ et $u\in K(S)$ une section globale qui
ne s'annule pas en $p$. Soit $p\in S_{1} \subset S$ un sous-sch\'ema de 
codimension $1$ qui soit lisse en $p$ et tel que $u(p) \notin T_{S_{1},p}$. Par construction, 
il existe un ouvert de Zariski $p\in S'_{1} \subset S_{1}$ qui est lisse et tel que 
pour tout $s\in S'_{1}(\mathbb{C})$ le noyau de l'application de Kodaira-Spencer
$$T_{S_{1},s} \longrightarrow H^{1}(X_{s},T_{X_{s}})$$
soit de dimension $r-1$, o\`u $r$ est le rang de $K$. 

Le champ de vecteurs
$u$ d\'efinit une action locale du groupe additif $\mathbb{C}$
$$f : U\times W \longrightarrow S^{an},$$
o\`u $U$ est un disque ouvert de $\mathbb{C}$ et $W$ est un voisinage
ouvert de $p$ dans $S^{an}$. 
Notons $W_{1}=W\cap S_{1}^{an}$, et consid\'erons le morphisme induit
$$g : U \times W_{1} \subset U\times W \longrightarrow S^{an}.$$
Ce morphisme est \'etale en $p$, car $u(p)$ et $ T_{S^{an}_{1},p}$ engendrent
$T_{S^{an},p}$. Son image contient donc un voisinage de $p$,  $V_{p} \subset S^{an}$.

Soit $s=g(x,t)$ un point de $V_{p}$, avec $(x,t)\in U\times W_{1}$. Par construction, la famille analytique
$$X^{an}\times_{S^{an}}U\times \{t\} \longrightarrow U$$
est telle que le morphisme de Kodaira-Spencer
$$T_{U,u} \longrightarrow H^{1}(X_{(u,t)},T_{X_{(u,t)}})$$
est nul pour tout $u\in U$. Comme de plus le rang de 
$H^{1}(X_{(u,t)},T_{X_{(u,t)}})$ est ind\'ependant de $u$, 
cette famille
est analytiquement isotriviale sur $U$ (voir \cite{kosp}), et comme $U$  est connexe
$X_{(u,t)}^{an}$  est ainsi isomorphe \`a $X_{(0,t)}^{an}$. Par GAGA on trouve donc
que $X_{(u,t)}\simeq X_{s}$ est isomorphe \`a $X_{(0,t)}\simeq X_{t}$. 
Nous avons donc montr\'e l'existence du sous-sch\'ema $S_{1}\subset S$, 
d'un ouvert de Zariski $S_{1}' \subset S_{1}$ lisse et contenant $p$,  
et d'un voisinage  $V_{p}$ de $p$ dans $S^{an}$, qui v\'erifient les deux
conditions suivantes.

\begin{enumerate}

\item Pour tout $s\in V_{p}$ il existe 
$t\in V_{p}\cap S_{1}^{an}$ tel que 
$X_{s}$ soit isomorphe \`a $X_{t}$.

\item Pour tout $s\in S'_{1}(\mathbb{C})$ le noyau de  l'application de Kodaira-Spencer
$$T_{S_{1},s} \longrightarrow H^{1}(X_{s},T_{X_{s}})$$
est de dimension $r-1$. 

\end{enumerate}

Quitte \`a restreindre $V_{p}$ on pourra aussi supposer que
$V_{p}\cap S_{1}^{an} \subset (S_{1}')^{an}$.  \\

Restreingons maintenant la famille $X$ sur l'ouvert $S_{1}'$ de $S_{1}$. En renouvelant 
plusieurs fois la m\^eme construction on d\'emontre que pour tout $p\in S^{an}$, il existe
un sous-sch\'ema $p\in S_{p} \subset S$, un ouvert de Zariski
$S_{p}' \subset S_{p}$ lisse et contenant $p$,  et 
un ouvert $p\in V_{p} \subset S^{an}$ avec $V_{p}\cap S_{1}^{an} \subset (S_{1}')^{an}$, et
qui v\'erifient les deux conditions suivantes.

\begin{enumerate}

\item Pour tout $s\in V_{p}$ il existe 
$t\in (S_{p}')^{an}$ tel que 
$X_{s}$ soit isomorphe \`a $X_{t}$.

\item Pour tout $s\in S'_{p}$ l'application de Kodaira-Spencer
$$T_{S'_{p},s} \longrightarrow H^{1}(X_{s},T_{X_{s}})$$
est injective. 

\end{enumerate}

Pour terminer, l'espace topologique $S^{an}$ se plonge comme sous-espace
ferm\'e dans $\mathbb{C}^{n}$, et est donc une r\'eunion d\'enombrable
de sous-espaces compacts (e.g. les traces des boules ferm\'ees de rayon
entier dans $\mathbb{C}^{n}$). Ceci implique qu'il existe un ensemble
d\'enombrable de points $\{p_{i}\}_{i\in I}$ de $S^{an}$ tel que les
ouverts $V_{p_{i}}$ recouvrent $S^{an}$. L'ensemble de familles
$$\{X\times_{S}S'_{p_{i}} \longrightarrow S_{p_{i}}'\}_{i\in I}$$
v\'erifie alors la conclusion du lemme. 
\hfill $\Box$ \\

\begin{lem}\label{l5}
Il existe un ensemble dominant de familles lisses, projectives et connexes
$\{X^{\alpha} \longrightarrow S^{\alpha}\}_{\alpha \in A}$ avec $A$ un ensemble
d\'enombrable, et tel que chaque famille $X^{\alpha} \longrightarrow S^{\alpha}$
soit miniverselle.
\end{lem}

\textit{Preuve du lemme \ref{l5}:} C'est une cons\'equence imm\'ediate
des lemmes \ref{l3} et \ref{l4}. \hfill $\Box$ \\

Le th\'eor\`eme \ref{t2} d\'ecoule maintenant du lemme \ref{l5} et du 
corollaire \ref{c2}. 
\hfill $\Box$ \\

\section{Remarques et compl\'ements}

Nous terminerons ce travail par quelques remarques et compl\'ements concernant les
r\'esultats pr\'esent\'es ainsi que les m\'ethodes utilis\'ees pour les d\'emontrer. \\

\begin{enumerate}

\item Comme nous l'avons mentionn\'e lors de l'introduction, la conjecture de 
Kawamata \ref{conj} implique le th\'eor\`eme \ref{t1}. En effet, si 
$X \longrightarrow S$ est une famille miniverselle de sch\'emas propres, lisses et connexes,  
et si $X_{0}$ est un sch\'ema propre et lisse sur un corps $K$, 
alors le sous-ensemble $S(X_{0})\subset S(K)$, form\'e des 
points $s$ tels que $X_{s}$ soit isomorphe \`a $X_{0}$, est (au plus) d\'enombrable (cela montre bien
que le conjecture \ref{conj} implique le th\'eor\`eme \ref{t1}). 
Ce fait se d\'emontre de la m\^eme fa\c{c}on que le th\'eor\`eme
\ref{t1}. On commence par consid\'erer le morphisme $S \longrightarrow \mathcal{VAR}^{c}_{smpr}$, 
et on forme le carr\'e homotopiquement cart\'esien suivant
$$\xymatrix{
S \ar[r] & \mathcal{VAR}^{c}_{smpr} \\
F \ar[r] \ar[u] & \mathcal{X}_{0}, \ar[u]}$$
o\`u $\mathcal{X}_{0}$ est le sous-champ plein de $\mathcal{VAR}^{c}_{smpr}$
qui est l'image du morphisme $X_{0} : Spec\, K \longrightarrow \mathcal{VAR}^{c}_{smpr}$. 
L'existence du sch\'ema de Hilbert implique alors que $F$ est un espace alg\'ebrique
localement de type fini sur $K$. De plus, les sch\'emas de Hilbert \'etant des r\'eunions
d\'enombrables de sous-espaces de type fini, on en d\'eduit que $F$ est
recouvert par un nombre d\'enombrable de sch\'emas affines. Enfin, la famille
$X\longrightarrow S$ \'etant miniverselle on voit que $F \longrightarrow Spec\, K$ est 
non-ramifi\'e, et donc de la forme $\coprod_{A}Spec\, K$ pour $A$ un ensemble
d\'enombrable. Ainsi, 
$S(X_{0})$, qui est l'image de $F(K) \longrightarrow S(K)$, est (au plus) d\'enombrable.

\item Le pr\'echamp $\mathcal{DGC}_{sat}^{c}$ poss\`ede un champ 
associ\'e $\widetilde{\mathcal{DGC}_{sat}^{c}}$. La proposition 
\ref{p1} montre que la diagonale de ce champ est repr\'esentable. Cependant, 
il semble que le champ $\widetilde{\mathcal{DGC}_{sat}^{c}}$ ne soit pas un champ
alg\'ebrique (ainsi, la r\'eponse \`a la question \cite[Q. 4.5]{to3} semble n\'egative). 
En effet, l'existence de d\'eformations formelles de sch\'emas propres et lisses 
qui ne sont pas alg\'ebrisables implique probablement l'existence
de d\'eformations formelles de dg-cat\'egories satur\'ees qui ne sont
pas alg\'ebrisables (\`a traver la construction $L_{pe}(-)$). Cela semble li\'e au fait qu'un sch\'ema 
(propre et lisse) formel
dont la cat\'egorie d\'eriv\'ee est satur\'ee est en fait alg\'ebrisable, bien que nous
ne savons pas si ce fait est toujours vrai (voir cependant \cite{bv} pour
un \'enonc\'e du m\^eme ordre dans le cadre analytique).

Le champ $\widetilde{\mathcal{DGC}_{sat}^{c}}$ se comporte en r\'ealit\'e de fa\c{c}on tout
\`a fait similaire \`a $\mathcal{VAR}^{c}_{smpr}$. Il est bien connu qu'une fa\c{c}on de corriger 
le fait que $\mathcal{VAR}^{c}_{smpr}$ ne soit pas un champ alg\'ebrique est 
de consid\'erer un champ de vari\'et\'es polaris\'ees.  De m\^eme, il semble que pour
obtenir un champ alg\'ebrique il faille consid\'erer un champ de dg-cat\'egories satur\'ees
munies de structures additionelles de type polarisation. Une piste possible serait de
regarder le champ des dg-cat\'egories satur\'ees munies de structure de stabilit\'e
\`a la Bridgeland. D\'efinir une version de $\widetilde{\mathcal{DGC}_{sat}^{c}}$ qui
soit un champ alg\'ebrique nous semble dans tous les cas une question particuli\`erement 
int\'eressante.

\item Pour terminer, signalons aussi l'existence d'un ''champ'' classifiant les 
dg-cat\'egories satur\'ees qui ne sont pas forc\'ement connexes. 
L'existence de groupes de cohomologie de Hochschild n\'egatifs non triviaux implique
que ce ''champ'' est en r\'ealit\'e un champ sup\'erieur (voir \cite[\S 4.3 (6)]{to3}, comme pr\'ec\'edemment, ce champ
sup\'erieur n'est probablement pas alg\'ebrique). Ce champ poss\`ede aussi une version 
''d\'eriv\'ee'', i.e. une version qui soit un $D^{-}$-champ au sens de \cite{hagII}. On s'attend \`a ce que le 
complexe tangent de ce champ pris en une dg-cat\'egorie $T$ soit $HH(T)[2]$, le complexe
de Hochschild d\'ecal\'e de $2$. En particulier, l'espace tangent en $T$ doit \^etre
$HH^{2}(T)$, ce qui semble \^etre confirm\'e par le travail en cours \cite{ke2}. 

L'identification de l'espace tangent de $\mathcal{DGC}_{sat}^{c}$ en $T$ avec 
$HH^{2}(T)$ permet aussi de donner une autre preuve de la propri\'et\'e de Torelli 
infinit\'esimal \ref{p2}, tout au moins lorsque l'on se place en caract\'eristique nulle. 
En effet, pour $T=L_{pe}(X)$, avec $X \longrightarrow Spec\, K$ un sch\'ema propre et lisse, 
et $K$ de caract\'eristique nulle, on a  (voir \cite{ye})
$$T^{0}_{T}\mathcal{DGC}_{sat}^{c}\simeq HH^{2}(T)\simeq H^{2}(X,\mathcal{O}_{X})\oplus 
H^{1}(X,T_{X}) \oplus H^{0}(X,\wedge^{2}T_{X}).$$
On voit ainsi que le morphisme
$$T^{0}\mathcal{VAR}^{c}_{smpr}\simeq H^{1}(X,T_{X}) \longrightarrow 
T^{0}_{T}\mathcal{DGC}_{sat}^{c}\simeq H^{2}(X,\mathcal{O}_{X})\oplus 
H^{1}(X,T_{X}) \oplus H^{0}(X,\wedge^{2}T_{X})$$
est un facteur direct. 

\end{enumerate}

\end{document}